\RequirePackage{fix-cm}
\documentclass[smallextended]{svjour3}       
\smartqed  
\usepackage{graphicx}
\usepackage[dvips]{color}
\usepackage{amsmath}
\usepackage{amssymb}
\usepackage{epsfig}
\usepackage{mathrsfs}
\usepackage{amsfonts}
\def\pn{\par\smallskip\noindent}
\def\proof{\pn {Proof.} }
\def\endproof{\hfill \quad{$\Box$}\smallskip}
\begin{document}

\title{A new semidefinite relaxation for $\ell_{1}$-constrained quadratic
optimization and extensions
\thanks{This research was supported by National
Natural Science Foundation of China under grants 11001006  and
91130019/A011702, by the fundamental research funds for the central universities
under grant
YWF-13-A01, and by the fund of State Key Laboratory of
Software Development Environment under grant SKLSDE-2013ZX-13.}
}

\titlerunning{A New SDP Bound for $\ell_{1}$-constrained QP}

\author{Yong Xia \and Yu-Jun Gong \and Sheng-Nan Han
}


\institute{Y. Xia \at
              State Key Laboratory of Software Development
              Environment, LMIB of the Ministry of Education, School of
Mathematics and System Sciences, Beihang University,  Beijing
100191,  P. R. China
              \email{dearyxia@gmail.com }
              \and
              Y.-J. Gong \and S.-N. Han  \at
School of Mathematics and System Sciences, Beihang University,  Beijing
100191,  P. R. China
              \email{gyjgongzuo@163.com, ~hsn20$\_$12@163.com}  
}

\date{Received: date / Accepted: date}

\maketitle

\begin{abstract}
In this paper, by improving the variable-splitting approach,
we propose a new semidefinite programming (SDP)
relaxation for the nonconvex quadratic optimization problem over the $\ell_1$
unit ball (QPL1).
It dominates the state-of-the-art SDP-based bound for (QPL1). As extensions,
we apply the new approach to the relaxation problem of
the sparse principal component analysis and
the nonconvex quadratic optimization problem over the $\ell_p$ ($1< p<2$)
unit ball and then show the dominance of the new relaxation.

\end{abstract}

\section{Introduction}
\label{intro}
We consider the quadratic optimization problem over the $\ell_1$ unit ball
\begin{eqnarray}
{\rm (QPL1(Q))}~ &\max& x^TQx \nonumber\\
 &{\rm s.t.}& \|x\|_{1}\leq 1, \nonumber
 \end{eqnarray}
which is known as an $\ell_1$-norm
trust-region subproblem in nonlinear programming \cite{C00} and $\ell_1$
Grothendieck problem in combinatorial optimization \cite{K12,K10}.
Applications of (QPL1(Q)) can be also found in compressed sensing
where $\|x\|_1$ is introduced to approximate $\|x\|_0$,
the number of nonzero elements of $x$.

If $Q$ is negative or positive semidefinite, (QPL1(Q)) is trivial to solve,
see \cite{Pin06}. Generally, (QPL1(Q)) is NP-hard, even when the
off-diagonal elements of $Q$ are all nonnegative, see \cite{H14}.
In the same paper,
Hsia showed that (QPL1(Q)) admits an exact nonconvex semidefinite
programming (SDP) relaxation, which was firstly proposed as an open problem
by Pinar and Teboulle
\cite{Pin06}.

Very recently, different SDP relaxations for (QPL1(Q)) have been studied
in \cite{Xia}.
The tightest one is the following doubly nonnegative (DNN) relaxation
due to Bomze et al. \cite{Bomze}:
\begin{eqnarray}
{\rm (DNN_{L1}(\tilde{Q}))}~ &\max& \tilde{Q}\bullet Y \nonumber\\
 &{\rm s.t.}& e^TYe=1, \nonumber\\
&&Y\ge 0,~Y\succeq 0,~Y\in {\cal S}^{2n}\nonumber
\end{eqnarray}
where $e$ is the vector with all elements equal to $1$, ${\cal S}^{2n}$
is the set of $2n\times 2n$ symmetric matrices, $Y\geq
0$  means that $Y$ is componentwise nonnegative, $Y\succeq 0$
stands for that $Y$ is positive semidefinite,   $A\bullet B={\rm
trace}(AB^T)=\sum_{i,j=1}^na_{ij}b_{ij}$ is the standard
inner product of $A$ and $B$, and
\[
\widetilde{Q}=\left[\begin{array}{cc}Q&-Q\\-Q&Q\end{array}\right].
\]
Notice that the set of extreme points of $\{x:~\|x\|_1\le 1\} $ is
$\{e_1,-e_1,\cdots,e_n,-e_n\}$,
where $e_i$ is the $i$-th column of the identity matrix $I$. Define
\[
A=[e_1,\cdots,e_n,-e_1,\cdots,-e_n]=[I~-I] \in \Re^{n\times 2n}.
\]
Then we have
\begin{equation}
 \{x\in \Re^n:~\|x\|_1\le 1\}=\{Ay:~e^Ty=1,y\ge 0,y\in  \Re^{2n}\}. \label{x-y}
\end{equation}
Consequently, (QPL1(Q)) can be equivalently transformed to
the following standard quadratic program (QPS) \cite{B02}:
\begin{eqnarray}
{\rm (QPS)}~ &\max_{y\in\Re^{2n}}& y^T\tilde{Q}y \nonumber\\
 &{\rm s.t.}& e^Ty=1,~y\ge 0. \nonumber
 \end{eqnarray}
Now we can see that ${\rm (DNN_{L1}(\tilde{Q}))}$ exactly corresponds to the well-known
 doubly nonnegative
relaxation of (QPS) \cite{Bomze}.
Moreover, as mentioned in \cite{Xia}, ${\rm (DNN_{L1}(\tilde{Q}))}$ can be
also derived by applying the
lifting procedure \cite{L90} to the following homogeneous reformulation of (QPS):
\begin{eqnarray}
&\max_{y\in\Re^{2n}}& y^T\tilde{Q}y \nonumber\\
 &{\rm s.t.}& y^Tee^Ty=1, \nonumber\\
 &&y_{i}y_{j}\ge 0,~i,j=1,\ldots,2n.\nonumber
 \end{eqnarray}

 A natural extension of  (QPL1(Q)) is
\begin{eqnarray}
({\rm QPL2L1(Q)})~ &\max & x^TQx \nonumber\\
&{\rm s.t.}&\|x\|_2=1, \nonumber \\
&&\|x\|_1^2\le k. \label{L1:1}
\end{eqnarray}
It is a relaxation of the sparse principal component analysis (SPCA)
problem \cite{L01} obtained by replacing the original constraint $\|x\|_0\le k$ with
(\ref{L1:1}) due to the following fact:
\[
\|x\|_1^2\le  \|x\|_0 \|x\|^2_2\le k.
\]
A well-known SDP relaxation  for (QPL2L1(Q)) is due to
d'Aspremont et al.  \cite{d07}:
\begin{eqnarray}
{\rm (SDP_X)}~&\max& Q\bullet X \nonumber\\
&{\rm s.t.}&  {\rm trace}(X)=1,\nonumber\\
&&e^T|X|e\le k, \nonumber\\
&&X\succeq 0,~X\in {\cal S}^{n}. \nonumber
\end{eqnarray}
Recently, Xia \cite{Xia} extended the doubly nonnegative relaxation approach
from (QPL1(Q)) to (QPL2L1(Q)) and obtained the following SDP relaxation:
\begin{eqnarray}
({\rm DNN_{\rm L2L1}(\widetilde{Q})})~ &\max & k\cdot\widetilde{Q}\bullet Y \nonumber\\
&{\rm s.t.}&k\cdot{\rm trace}( A^TAY)=1,  \nonumber\\
&&   e^TYe=1,  \nonumber\\
&&Y\ge 0,~Y\succeq 0,~Y\in {\cal S}^{2n}.\nonumber
\end{eqnarray}
It was proved in \cite{Xia} that $v({\rm DNN_{\rm L2L1}(\widetilde{Q})})=v{\rm (SDP_X)}$, where $v(\cdot)$ denote the optimal value of problem $(\cdot)$. Unfortunately,  this equivalence result is incorrect though it is true that  $v({\rm DNN_{\rm L2L1}(\widetilde{Q})})\le v{\rm (SDP_X)}$. A first counterexample will be given in this paper (see Example 2 below) to show it is possible $v({\rm DNN_{\rm L2L1}(\widetilde{Q})})<v{\rm (SDP_X)}$.

The other  extension of  (QPL1(Q)) is
\begin{eqnarray}
({\rm QPLp(Q)})~ &\max& x^TQx \nonumber\\
 &{\rm s.t.}& \|x\|_{p}\leq 1,  \nonumber
 \end{eqnarray}
where  $\|x\|_p = \left(\sum^n_{i=1}|x_{i}|^p\right)^{\frac{1}{p}}$ and $1< p< 2$.
(QPLp) is known as a special case of the $\ell_p$ Grothendieck problem
if the diagonal entries of $Q$ vanish.
According to the survey \cite{K12},
there is no approximation and hardness results for the $\ell_p$ Grothendieck
 problem with $1<p<2$.
Though $({\rm QPLp(Q)})$ has an exact nonconvex SDP relaxation similar to
that of $({\rm QPL1(Q)})$, the computational complexity of
$({\rm QPLp(Q)})$ is still unknown \cite{H14}.

Since the $\ell_p$ unit balls ($1<p<2$) are included in the $\ell_2$ unit ball,
a trivial bound for $({\rm QPLp(Q)})$ is
\begin{equation}\label{B2}
B_2(Q):=\max_{\|x\|_2\le 1} x^TQx= \max\left\{\lambda_{\max}(Q),0\right\},
\end{equation}
where $\lambda_{\max}(Q)$ is the largest eigenvalue of $Q$.

As mentioned by Nesterov in the SDP Handbook \cite{Nes00},
no practical SDP bounds of $({\rm QPLp(Q)})$ are in sight for $1<p<2$.
Recently, Bomze \cite{Bomze} used the H\"{o}lder inequality
\begin{equation}
\|x\|_1 \leq \|x\|_p\|e\|_{\frac{p}{p-1}} = n^{\frac{p-1}{p}}\|x\|_p\label{hol}
\end{equation}
to propose the following SDP bound
\begin{equation}
B_1(Q):=n^{\frac{2(p-1)}{p}}\cdot v{\rm (DNN_{L1}(\tilde{Q}))}. \label{Bom}
\end{equation}
In general, ${\rm B_1(Q)}$ dominates ${\rm B_2(Q)}$ when $p$ close to $1$,
though lacking a proof.

In this paper,
 based on a new variable-splitting reformulation for
 the $\ell_1$-constrained set,
we establish a new SDP relaxation for (QPL1(Q)), which is proved to dominate
${\rm (DNN_{L1}(\tilde{Q}))}$.
We use a small example to show the improvement could be strict.
Then we extend the new approach to (QPL2L1(Q)) and obtain two new SDP relaxations.
We cannot prove the first new SDP bound dominates
${\rm (DNN_{L2L1}(\tilde{Q}))}$, though it was demonstrated by examples.
However, under a mild assumption, the second new SDP bound dominates
${\rm (DNN_{L2L1}(\tilde{Q}))}$.
Finally, motivated by the model (QPL2L1(Q)),
we establish a new SDP bound for (QPLp(Q)) and show it is in general
tighter than $\min\{{\rm B_2(Q)}, {\rm B_1(Q)}\}$.

The paper is organized as follows. In Section 1, we propose
a new variable-splitting reformulation for the $\ell_1$-constrained set
and then a new SDP relaxation for (QPL1(Q)). We show it improves
the state-of-the-art SDP-based bound.
In Section 2, we extend the new SDP approach to (QPL2L1(Q)) and study the obtained
two new
SDP relaxations. In Section 3,
we establish a new SDP relaxation for (QPLp(Q)), which
improves the existing upper bounds.
Conclusions are made in Section 4.

\section{A New SDP Relaxation  for (QPL1(Q))}
In this section,we establish a new SDP relaxation for
(QPL1(Q)) based on a new variable-splitting reformulation for
 the $\ell_1$-constrained set.

For any $x\in \Re^n$, let
 \begin{eqnarray*}
&&y_{i}=\max\{x_i,0\},~i=1,\ldots,n,\\
&&y_{n+i}=-\min\{x_i,0\},~i=1,\ldots,n.
 \end{eqnarray*}
Then we have
 \begin{eqnarray}
&&x_{i}=y_{i}-y_{i+n},~i=1,\ldots,n,\label{x-y:1}\\
&&|x_{i}|=y_{i}+y_{i+n},~i=1,\ldots,n,\label{x-y:2}\\
&&y_{i} y_{i+n}=0,~i=1,\ldots,n, \label{x-y:3}\\
&&y_{i}\ge 0,~i=1,\ldots,2n. \label{x-y:4}
 \end{eqnarray}
Now we obtain a new variable-splitting
reformulation of the $\ell_1$-constrained set:
\begin{equation*}
 \{x:~\|x\|_1\le 1\}=\{Ay:~ e^Ty\le 1,y\ge 0,y\in \Re^{2n},y_{i} y_{i+n}=0,
 i=1,\ldots,n\}.
\end{equation*}
It follows that
\begin{eqnarray}
v({\rm QPL1(Q)})&=&\max_{y\in \Re^{2n}} y^T\tilde{Q}y \nonumber\\
 &&~~{\rm s.t.}~ e^Ty\le1,y\ge 0,\nonumber\\
  &&~~~~~~~y_{i} y_{i+n}=0,~i=1,\ldots,n, \nonumber\\
 &=&\max_{y\in \Re^{2n}} y^T\tilde{Q}y \nonumber\\
 &&~~{\rm s.t.}~ e^Tyy^Te\le1, \nonumber\\
 &&~~~~~~~y_{i} y_{i+n}=0,~i=1,\ldots,n, \nonumber\\
  &&~~~~~~~y_{i} y_{j}\ge0,~i,j=1,\ldots,2n. \nonumber
 \end{eqnarray}
Applying the
lifting procedure \cite{L90},
we obtain the following new doubly nonnegative relaxation of (QPL1(Q))
 \begin{eqnarray*}
{\rm (DNN^{\rm new}_{L1}(\tilde{Q}))}~ &\max& \tilde{Q}\bullet Y \nonumber\\
 &{\rm s.t.}&  e^TYe\leq1,  \\
  &&Y_{i,n+i} = 0,~ i=1,\ldots,n, \\
&&Y\ge0,~Y\succeq 0,~Y\in {\cal S}^{2n}.
 \end{eqnarray*}
We first compare the qualities of  $v({\rm DNN_{L1}})$
and $v({\rm DNN^{\rm new}_{L1}})$.
\begin{theorem}\label{thm:1}
$v({\rm DNN_{L1}(\tilde{Q})}) \geq v({\rm DNN^{\rm new}_{L1}(\tilde{Q})})
\ge v({\rm QPL1(Q)})$.
\end{theorem}
\proof
According to the definitions, we have $v({\rm DNN_{L1}}(\tilde{Q}))
\ge v({\rm QPL1(Q)})$ and
$v({\rm DNN^{\rm new}_{L1}(\tilde{Q})}) \ge v({\rm QPL1(Q)})$.
It is sufficient to prove the first inequality.

Since $Y=0_{2n\times 2n}$ is a feasible solution of
$({\rm DNN^{\rm new}_{L1}(\tilde{Q})})$, we have
\[
v({\rm DNN^{\rm new}_{L1}(\tilde{Q})})\ge 0.
\]
Suppose $Q\preceq 0$.
Let $Y^*$ be an optimal solution of (${\rm DNN^{\rm new}_{L1}(\tilde{Q})}$).
  Since $Y^*\succeq 0$, we have $AY^*A^T\succeq 0$ and therefore
\[
v({\rm DNN^{\rm new}_{L1}(\tilde{Q})})=\tilde{Q}\bullet Y^*
={\rm trace}((A^TQA)Y^*)={\rm trace}(Q(AY^*A^T))\le 0.
\]
Consequently, $v({\rm DNN^{\rm new}_{L1}(\tilde{Q})})=0$.
Similarly, we can show $v({\rm DNN_{L1}(\tilde{Q})})=0$.

Now we assume $Q\not\preceq 0$. There is a vector $v$
such that $\|v\|_1\le 1$ and $v^TQv>0$.
That is, $v({\rm QPL1(Q)})>0$. It follows that
$v({\rm DNN^{\rm new}_{L1}(\tilde{Q})}) >0$.
Let $Y^*$ be an optimal solution of (${\rm DNN^{\rm new}_{L1}(\tilde{Q})}$).
Then $Y^*\neq 0_{2n\times 2n}$. Moreover, since $Y^*\ge 0$,
we have $e^TY^*e>0$. We conclude that
\begin{equation}
e^TY^*e=1.\label{eq1}
\end{equation}
If this is not true, then $0<e^TY^*e<1$. Define \[
\widetilde{Y}=\frac{1}{e^TY^*e}Y^*.
\]
It is trivial to see that $\widetilde{Y}$ is also feasible to
(${\rm DNN^{\rm new}_{L1}(\tilde{Q})}$). Moreover, we have
\[
\tilde{Q}\bullet \widetilde{Y}=\frac{1}{e^TY^*e}\tilde{Q}\bullet Y^*
>\tilde{Q}\bullet Y^*,
\]
which contradicts the fact that $Y^*$ is a maximizer of
(${\rm DNN^{\rm new}_{L1}(\tilde{Q})}$).
According to the equality (\ref{eq1}),
$Y^*$ is also a feasible solution of (${\rm DNN_{L1}(\tilde{Q})}$).
Consequently, $v({\rm DNN_{L1}(\tilde{Q})}) \geq
v({\rm DNN^{\rm new}_{L1}(\tilde{Q})})$. The proof is complete.
\endproof

The following small example illustrates that
$v({\rm DNN^{\rm new}_{L1}(\tilde{Q})})$ could strictly improve
$v({\rm DNN_{L1}(\tilde{Q})})$.
\begin{example}\label{exam1}
Consider the following instance of dimension $n=6$
\[
Q=\left[\begin{array}{cccccc}
   -11&   -11 &   -7&   -10&    -8&    -2\\
   -11&    -5 &  -10&    -9&   -10&    -7\\
    -7&   -10 &  -10&    -3&    -6&    -8\\
   -10 &   -9 &   -3&    -8&    -9&   -10\\
    -8 &  -10 &   -6&    -9&    -8&    -7\\
    -2 &   -7 &   -8&   -10&    -7&    -6\end{array}\right]
    \]
We modeled this instance by  CVX 1.2  (\cite{CVX}) and solved it by
SEDUMI (\cite{St}) within CVX. Then we obtained that
\[v({\rm DNN_{L1}(\tilde{Q})})\thickapprox2.0487,~
v({\rm DNN^{\rm new}_{L1}(\tilde{Q})})\thickapprox2.0186.\]
\end{example}

Finally, we show that there are some cases for which
$({\rm DNN^{\rm new}_{L1}(\tilde{Q})})$ has no improvement. This ``negative''
result is also interesting in the sense that
in case we solve $({\rm DNN_{L1}(\tilde{Q})})$,
we can fix $Y_{i,n+i}$ ($i=1,\ldots,n$)
at zeros in advance.
\begin{theorem}\label{thm:2}
Suppose ${\rm diag}(Q)\ge 0$.
$v({\rm DNN_{L1}(\tilde{Q})}) = v({\rm DNN^{\rm new}_{L1}(\tilde{Q})})$.
\end{theorem}
\proof
Let $Y^*$ be an optimal solution of $({\rm DNN_{L1}(\tilde{Q})})$.
Suppose there is an index $k\in\{1,\ldots,n\}$ such that
$Y^*_{k,n+k}>0$. Let $\delta_k=Y^*_{k,n+k}$ and define a symmetric matrix
$Z\in {\cal S}^{2n}$ where
\[
Z_{kk}=Z_{n+k,n+k}=\delta_k, ~Z_{k,n+k}=Z_{n+k,k}=-\delta_k
\]
and all other elements are zeros. Then
\[
Z\succeq 0,~\tilde{Q}\bullet Z=2(Q_{kk}+Q_{n+k,n+k})\delta_k \ge 0.
\]
It follows that
\[
Y^*+Z\succeq 0,~ Y^*+Z\ge 0,~ (Y^*+Z)_{k,n+k}=0,~
\tilde{Q}\bullet (Y^*+Z)\ge \tilde{Q}\bullet Y^*.
\]
Then, $Y^*+Z$ is also
an optimal solution of $({\rm DNN_{L1}(\tilde{Q})})$.
Repeat the above procedure until we obtain an
optimal solution of $({\rm DNN_{L1}(\tilde{Q})})$,
denoted by $\widetilde{Y}^*$,
satisfying $\widetilde{Y}^*_{i,n+i}=0$ for $i=1,\ldots,n$.
Notice that $\widetilde{Y}^*$ is a feasible solution of
 $({\rm DNN^{\rm new}_{L1}})$. Therefore, we have
 $v({\rm DNN_{L1}(\tilde{Q})}) \le v({\rm DNN^{\rm new}_{L1}(\tilde{Q})})$.
Combining this inequality with Theorem \ref{thm:1}, we can complete the proof.
\endproof

\section{New SDP Relaxations for (QPL2L1(Q))}

In this section, we extend the above new reformulation approach to (QPL2L1(Q))
 and obtain two new semidefinite programming relaxations.

Similar to the
reformulation (\ref{x-y:1})-(\ref{x-y:4}), we have
 \begin{eqnarray}
&&x_{i}=\sqrt{k}(y_{i}-y_{n+i}),~i=1,\ldots,n,\label{xy:1}\\
&&|x_{i}|=\sqrt{k}(y_{i}+y_{n+i}),~i=1,\ldots,n,\label{xy:2}\\
&&y_{i} y_{n+i}=0,~i=1,\ldots,n, \label{xy:3}\\
&&y_{i}\ge 0,~i=1,\ldots,2n. \label{xy:4}
 \end{eqnarray}
It follows that
\begin{eqnarray}
 &&\{x:~\|x\|_2=1,~\|x\|_1\le k\}\nonumber\\
 &=&\{\sqrt{k}Ay:~ky^TA^TAy=1, e^Ty\le 1,y\ge 0,y\in \Re^{2n},
 y_{i} y_{n+i}=0,i=1,\ldots,n\}. \nonumber
\end{eqnarray}
Introducing $Y=yy^T\succeq 0$, we obtain the following new SDP relaxation
for (QPL2L1(Q)):
 \begin{eqnarray}
 ({\rm DNN^{\rm new\le}_{L2L1}(\tilde{Q})}):~
 &\max & k\cdot \widetilde{Q}\bullet Y  \nonumber\\
&{\rm s.t.}&k\cdot {\rm trace}(A^TAY)=1, \nonumber \\
&&    e^TYe \leq 1,  \nonumber\\
&&Y_{i,n+i}=0, ~i=1,\ldots,n,\nonumber\\
&&Y\ge0,~Y\succeq 0,~Y\in {\cal S}^{2n}.\nonumber
\end{eqnarray}
According to the definition, we trivially have:
\begin{proposition}
$ v({\rm DNN^{\rm new\le}_{L2L1}(\tilde{Q})})\ge v({\rm QPL2L1(Q)})$.
\end{proposition}
\begin{proposition}\label{prop}
$\max\left\{v({\rm DNN_{L2L1}(\tilde{Q})}),
~v({\rm DNN^{\rm new\le}_{L2L1}(\tilde{Q})})\right\}\le \lambda_{\max}(Q)$.
\end{proposition}
\proof
Both $({\rm DNN_{L2L1}(\tilde{Q})})$ and
$({\rm DNN^{\rm new\le}_{L2L1}(\tilde{Q})})$
share the same relaxation:
 \begin{eqnarray}
 ({\rm R_Y})~ &\max & k\cdot \widetilde{Q}\bullet Y  \nonumber\\
&{\rm s.t.}&k\cdot {\rm trace}(A^TAY)=1, \nonumber \\
&&Y\succeq 0.\nonumber
\end{eqnarray}
Let $X=kAYA^T$. We have
\begin{eqnarray*}
&&k\cdot \widetilde{Q}\bullet Y=Q\bullet X,\\
&&k\cdot{\rm trace}(A^TAY)={\rm trace}(X),\\
&&Y\succeq 0 \Longrightarrow X\succeq 0.
\end{eqnarray*}
 Therefore, $({\rm R_Y})$ can be further relaxed to
 \begin{eqnarray}
 ({\rm R_X})~ &\max & Q\bullet X  \nonumber\\
&{\rm s.t.}&  {\rm trace}(X)=1, \nonumber \\
&& X\succeq 0.\nonumber
\end{eqnarray}
Let $Q=U\Sigma U^T$ be the eigenvalue decomposition of $Q$,
where $\Sigma={\rm Diag}(\sigma_1,\ldots,\sigma_n)$ and $U$ are
 column-orthogonal. Since
\begin{eqnarray}
&&{\rm trace}(X)={\rm trace}(U^TXU),\label{x:1}\\
&&X\succeq 0\Longrightarrow X_{ii}\ge 0,\label{x:2}\\
&&X\succeq 0 \Longleftrightarrow U^TXU\succeq 0,\label{x:3}
\end{eqnarray}
 we can further relax  $({\rm R_X})$ to the following linear
 programming problem:
 \begin{eqnarray}
 ({\rm LP})~ &\max & \sum_{i=1}^n \sigma_ix_i  \nonumber\\
&{\rm s.t.}&  \sum_{i=1}^nx_i=1, \nonumber \\
&& x_i\ge 0,~i=1,\ldots,n.\nonumber
\end{eqnarray}
Now it is trivial to verify that
\[
v({\rm LP})=\max\{\sigma_1,\ldots,\sigma_n\}=\lambda_{\max}(Q).
\]
The proof is complete.
\endproof
\begin{corollary}\label{cor}
Suppose $v({\rm QPL2L1(Q)})=\lambda_{\max}(Q)$, then we have
\[
v({\rm DNN_{L2L1}(\tilde{Q})})=v({\rm DNN^{\rm new\le}_{L2L1}(\tilde{Q})})
=v({\rm QPL2L1(Q)}).
\]
\end{corollary}

We are unable to prove $v({\rm DNN_{L2L1}(\tilde{Q})}) \geq
v({\rm DNN^{\rm new\le}_{L2L1}(\tilde{Q})})$, though we failed to
have found an example such that $v({\rm DNN_{L2L1}(\tilde{Q})}) <
v({\rm DNN^{\rm new\le}_{L2L1}(\tilde{Q})})$.
Moreover, the following example shows that it is possible
$v({\rm DNN_{L2L1}(\tilde{Q})}) > v({\rm DNN^{\rm new\le}_{L2L1}(\tilde{Q})})$. As a by-product, we observe 
$v({\rm DNN_{\rm L2L1}(\widetilde{Q})})<v{\rm (SDP_X)}$ from the example, which means that the result 
$v({\rm DNN_{\rm L2L1}(\widetilde{Q})})=v{\rm (SDP_X)}$ (Theorem 3.2 \cite{Xia}) is incorrect. Notice that it is true that $v({\rm DNN_{\rm L2L1}(\widetilde{Q})})\le v{\rm (SDP_X)}$.
 
\begin{example}\label{exam2}
Consider the same instance of Example \ref{exam1} and let $k=3$.
We modeled this instance by  CVX 1.2  (\cite{CVX}) and solved it by
SEDUMI (\cite{St}) within CVX. We obtained that
\[
v{\rm (SDP_X)}\thickapprox6.3104,~
v({\rm DNN_{L2L1}(\tilde{Q})})\thickapprox6.0964,~
v({\rm DNN^{\rm new\le}_{L2L1}(\tilde{Q})})
\thickapprox5.9962.
\]
\end{example}

Thus, in order to theoretically improve $v({\rm DNN_{L2L1}(\tilde{Q})})$,
we consider
 \begin{eqnarray}
 ({\rm DNN^{\rm new=}_{L2L1}(\tilde{Q})}) ~ &\max& k\cdot\tilde{Q}\bullet Y \nonumber\\
 &{\rm s.t.}& k \cdot {\rm trace}(Y)=1,  \nonumber\\
&&  e^TYe =  1, \nonumber\\
  &&Y_{i,n+i} = 0, ~i=1,\ldots,n, \nonumber\\
&&Y\ge0,~Y\succeq 0,~Y\in {\cal S}^{2n}.\nonumber
 \end{eqnarray}
It is trivial to see that
\[
v({\rm DNN_{L2L1}(\tilde{Q})}) \geq v({\rm DNN^{\rm new=}_{L2L1}(\tilde{Q})}).
\]
However, $v({\rm DNN^{\rm new=}_{L2L1}(\tilde{Q})})$ may be not an upper bound
of $({\rm QPL2L1(Q)})$, which is indicated by the following example.
\begin{example}\label{exam3}
Consider the same instance of Example \ref{exam1} and let $k=5$.
We modeled this instance by  CVX 1.2  (\cite{CVX}) and solved it by
SEDUMI (\cite{St}) within CVX. We obtained that
\[
v({\rm DNN^{\rm new=}_{L2L1}(\tilde{Q})})\thickapprox7.048
<v({\rm QPL2L1(Q)})=\lambda_{\max}(Q)=7.0857.
\]
\end{example}
So,
we have to identify when $v({\rm DNN^{\rm new=}_{L2L1}(\tilde{Q})})$
is an upper bound of $({\rm QPL2L1(Q)})$.
\begin{theorem}\label{thm:4}
Suppose
\begin{equation}
v({\rm QPL2L1(Q)})<\lambda_{\max}(Q), \label{as:2}
\end{equation}
we have
$v({\rm DNN^{\rm new=}_{L2L1}(\tilde{Q})})\ge v({\rm QPL2L1(Q)})$.
\end{theorem}
\proof
We first notice that the maximum eigenvalue problem
\[
{\rm(E)}~\max_{\|x\|_2=1}x^TQx = \lambda_{\max}(Q)
\]
is a homogeneous trust-region subproblem and hence
has no local-non-global maximizer \cite{M94}.
Therefore,
suppose there is an optimal solution of $({\rm QPL2L1(Q)})$, denoted by $x^*$,
satisfying $\|x\|_1^2< k$, then $x^*$ also globally solves (E), i.e.,
\[
v({\rm QPL2L1(Q)})=x^{*T}Qx^*=\lambda_{\max}(Q).
\]
Consequently,  the assumption (\ref{as:2}) implies that
\begin{eqnarray}
v({\rm QPL2L1(Q)})=&\max & x^TQx \nonumber\\
&{\rm s.t.}&\|x\|_2=1 \nonumber \\
&&\|x\|_1^2= k.\nonumber
\end{eqnarray}
Taking the transformation (\ref{xy:1})-(\ref{xy:4}) and then applying the
lifting approach \cite{L90},  we obtain the SDP relaxation
$({\rm DNN^{\rm new=}_{L2L1}(\tilde{Q})})$.
The proof is complete.
\endproof

\begin{remark}
The assumption (\ref{as:2})
is generally not easy to verify. 
However, when $Q$ has a unique maximum eigenvalue, (\ref{as:2}) holds if and only
if $\|v\|_1>\sqrt{k}$, where $v$ is the $\ell_2$-normalized eigenvector corresponding to the
maximum eigenvalue of $Q$.
Moreover,
according to Corollary \ref{cor} and Proposition \ref{prop},
the assumption (\ref{as:2})
can be replaced by the following easy-to-check sufficient condition
\[
v({\rm DNN_{L2L1}(\tilde{Q})})<\lambda_{\max}(Q).
\]
\end{remark}


\section{A New SDP Relaxation for (QPLp(Q)) ($1< p<2$)}

In this section,
we first propose a new SDP relaxation for (QPLp(Q)) and then
show it improves both $B_2(Q)$ (\ref{B2}) and $B_1(Q)$ (\ref{Bom}).

Motivated by the  H\"{o}lder inequality (\ref{hol}) and
 the model (QPL2L1(Q)), we obtain the following new
  relaxation for (QPLp(Q)):
\begin{eqnarray}
({\rm QPL2L1^{\le}(Q)})~ &\max & x^TQx \nonumber\\
&{\rm s.t.}&\|x\|_2\le 1 \nonumber\\
&&\|x\|_1^2\le n^{\frac{2(p-1)}{p}}.\nonumber
\end{eqnarray}

Taking the transformation (\ref{xy:1})-(\ref{xy:4}) and then applying the
lifting approach \cite{L90},  we obtain the following SDP relaxation for
${\rm (QPL2L1^{\le}(Q))}$, which is very similar to
$({\rm DNN^{\rm new\le}_{L2L1}(\tilde{Q})})$:
 \begin{eqnarray}
 ({\rm DNN_{Lp}(\widetilde{Q})})~ &\max& n^{\frac{2(p-1)}{p}}\cdot
  \tilde{Q}\bullet Y \nonumber\\
 &{\rm s.t.}& n^{\frac{2(p-1)}{p}}\cdot  {\rm trace}(Y)\leq 1  \nonumber\\
&&  e^TYe \le  1, \nonumber\\
  &&Y_{i,n+i} = 0, ~i=1,\ldots,n, \nonumber\\
&&Y\ge0,~Y\succeq 0,~Y\in {\cal S}^{2n}.\nonumber
 \end{eqnarray}
\begin{theorem}
\[
\min\{{\rm B_2(Q),B_1(Q)}\}\ge v({\rm DNN_{Lp}(\widetilde{Q})})
\ge v({\rm QPLp(Q)}).
\]
\end{theorem}
\proof
According to the definitions, the second inequality is trivial.
It is sufficient to prove the first inequality. We first show
$ {\rm B_2(Q)}  \ge v({\rm DNN_{Lp}(\widetilde{Q})})$.

Let $X=n^{\frac{2(p-1)}{p}}AYA^T$. Since
\begin{eqnarray*}
&&n^{\frac{2(p-1)}{p}}\cdot \tilde{Q}\bullet Y= Q\bullet X,\\
&&n^{\frac{2(p-1)}{p}}\cdot{\rm trace}(A^TAY)={\rm trace}(X),\\
&&Y\succeq 0 \Longrightarrow X\succeq 0,
\end{eqnarray*}
$({\rm DNN_{Lp}(\widetilde{Q})})$ has the following relaxation:
\begin{eqnarray}
 ({\rm R})~ &\max & Q\bullet X  \nonumber\\
&{\rm s.t.}&  {\rm trace}(X)\le 1, \nonumber \\
&& X\succeq 0.\nonumber
\end{eqnarray}
Let $Q=U\Sigma U^T$ be the eigenvalue decomposition of $Q$,
where $\Sigma={\rm Diag}(\sigma_1,\ldots,\sigma_n)$ and $U$ are
 column-orthogonal. According to (\ref{x:1})-(\ref{x:3}),
 we can further relax  $({\rm R})$ to the following linear
 programming problem:
 \begin{eqnarray}
 ({\rm LP})~ &\max & \sum_{i=1}^n \sigma_ix_i  \nonumber\\
&{\rm s.t.}&  \sum_{i=1}^nx_i\le 1, \nonumber \\
&& x_i\ge 0,~i=1,\ldots,n.\nonumber
\end{eqnarray}
It is not difficult  to verify that
\[
v({\rm LP})=\max\{0,\sigma_1,\ldots,\sigma_n\}=\max\{0,\lambda_{\max}(Q)\}=
{\rm B_2(Q)}.
\]
Now we prove $ {\rm B_1(Q)}  \ge v({\rm DNN_{Lp}(\widetilde{Q})})$. Notice that
 \begin{eqnarray}
 n^{-\frac{2(p-1)}{p}}\cdot v({\rm DNN_{Lp}(\widetilde{Q})})&\le &\max
  \tilde{Q}\bullet Y \nonumber\\
 &~~~~&{\rm s.t.}~  e^TYe \le  1, \nonumber\\
  &~~~~~~&~~~~~Y_{i,n+i} = 0, ~i=1,\ldots,n, \nonumber\\
&~~~~~~&~~~~~Y\ge0,~Y\succeq 0,~Y\in {\cal S}^{2n}\nonumber\\
&=&v{\rm (DNN^{\rm new}_{L1}(\tilde{Q}))}\nonumber\\
&\le& v{\rm (DNN_{L1}(\tilde{Q}))},\nonumber
 \end{eqnarray}
 where the last inequality follows from Theorem \ref{thm:1}.
  The proof is complete.
\endproof

We randomly generated a symmetric matrix $Q$ of order $n=10$
using the following Matlab scripts:
\begin{verbatim}
rand('state',0); Q = rand(n,n); Q = (Q+Q')/2;
\end{verbatim}
and then compared the qualities of the three upper bounds,
$v({\rm DNN_{Lp}(\widetilde{Q})})$, ${\rm B_1(Q)}$ and ${\rm B_2(Q)}$.
The results were plotted in Figure 1, where
the lower bound of QPLp(Q) is computed as follows.
Solve $({\rm DNN_{Lp}(\widetilde{Q})})$ and obtain the optimal solution $Y^*$.
Let $y, z$ be the unit eigenvectors corresponding to the maximum eigenvalues
 of $AY^*A^T$ and $Q$, respectively.
Then $\frac{1}{\|y\|_p}y$ and $\frac{1}{\|z\|_p}z$
are two feasible solutions of (QPLp(Q)) and
\[
\max\left\{\frac{y^TQy}{\|y\|_p^2},~\frac{z^TQz}{\|z\|_p^2}\right\}
\]
gives a lower bound of $v({\rm QPLp(Q)})$.
From Figure 1, we can see that for $1<p<2$,
though ${\rm B_2(Q)}$ and ${\rm B_1(Q)}$ cannot dominate each other,
both are strictly improved by
$v({\rm DNN_{Lp}(\widetilde{Q})})$.

\begin{figure}\label{fig}
\begin{center}
\includegraphics[width=0.8\textwidth]{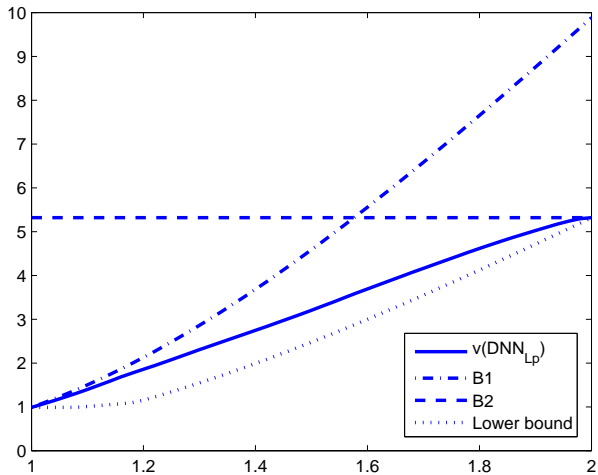}
\caption{Quality of the lower bound and the three upper bounds, $B_2(Q)$, $B_1(Q)$ and
$v({\rm DNN_{Lp}(\widetilde{Q})})$ in dependence of $p$.}
\end{center}
\end{figure}

\section{Conclusion}
The SDP relaxation has been known to generate high quality bounds for
nonconvex quadratic optimization problems. In this paper, based on
a new variable-splitting characterization of the $\ell_1$
unit ball, we establish a new semidefinite programming (SDP)
relaxation for the quadratic optimization problem over the $\ell_1$
unit ball (QPL1).
We show the new developed SDP bound dominates
the state-of-the-art SDP-based upper bound for (QPL1). There is an example to
show the improvement could be strict.
Then we extend the new reformulation approach to
the relaxation problem of
the sparse principal component analysis (QPL2L1) and obtain two SDP
formulations.
Examples demonstrate that the first SDP bound is in general tighter than
the DNN relaxation for (QPL2L1). But we are unable to prove it.
Under a mild assumption, the second SDP bound dominates the DNN relaxation.
Finally, we extend our approach to
the nonconvex quadratic optimization problem over the $\ell_p$ ($1< p<2$)
unit ball (QPLp) and show the new SDP bound dominates
two upper bounds in recent literature.
%



\end{document}